\input amstex
%Sekretariats-Stildatei
%
%AMSPPT.STY modifiziert und erweitert von C. Krattenthaler

\def\next{AMS-SEKR}\ifx\styname\next \endinput\fi
\catcode`\@=11
\def\styname{AMS-SEKR}
\def\styversion{2.0}
{\W@{}\W@{\styname.STY - Version \styversion}\W@{}}
\hyphenation{acad-e-my acad-e-mies af-ter-thought anom-aly anom-alies
an-ti-deriv-a-tive an-tin-o-my an-tin-o-mies apoth-e-o-ses apoth-e-o-sis
ap-pen-dix ar-che-typ-al as-sign-a-ble as-sist-ant-ship as-ymp-tot-ic
asyn-chro-nous at-trib-uted at-trib-ut-able bank-rupt bank-rupt-cy
bi-dif-fer-en-tial blue-print busier busiest cat-a-stroph-ic
cat-a-stroph-i-cally con-gress cross-hatched data-base de-fin-i-tive
de-riv-a-tive dis-trib-ute dri-ver dri-vers eco-nom-ics econ-o-mist
elit-ist equi-vari-ant ex-quis-ite ex-tra-or-di-nary flow-chart
for-mi-da-ble forth-right friv-o-lous ge-o-des-ic ge-o-det-ic geo-met-ric
griev-ance griev-ous griev-ous-ly hexa-dec-i-mal ho-lo-no-my ho-mo-thetic
ideals idio-syn-crasy in-fin-ite-ly in-fin-i-tes-i-mal ir-rev-o-ca-ble
key-stroke lam-en-ta-ble light-weight mal-a-prop-ism man-u-script
mar-gin-al meta-bol-ic me-tab-o-lism meta-lan-guage me-trop-o-lis
met-ro-pol-i-tan mi-nut-est mol-e-cule mono-chrome mono-pole mo-nop-oly
mono-spline mo-not-o-nous mul-ti-fac-eted mul-ti-plic-able non-euclid-ean
non-iso-mor-phic non-smooth par-a-digm par-a-bol-ic pa-rab-o-loid
pa-ram-e-trize para-mount pen-ta-gon phe-nom-e-non post-script pre-am-ble
pro-ce-dur-al pro-hib-i-tive pro-hib-i-tive-ly pseu-do-dif-fer-en-tial
pseu-do-fi-nite pseu-do-nym qua-drat-ics quad-ra-ture qua-si-smooth
qua-si-sta-tion-ary qua-si-tri-an-gu-lar quin-tes-sence quin-tes-sen-tial
re-arrange-ment rec-tan-gle ret-ri-bu-tion retro-fit retro-fit-ted
right-eous right-eous-ness ro-bot ro-bot-ics sched-ul-ing se-mes-ter
semi-def-i-nite semi-ho-mo-thet-ic set-up se-vere-ly side-step sov-er-eign
spe-cious spher-oid spher-oid-al star-tling star-tling-ly
sta-tis-tics sto-chas-tic straight-est strange-ness strat-a-gem strong-hold
sum-ma-ble symp-to-matic syn-chro-nous topo-graph-i-cal tra-vers-a-ble
tra-ver-sal tra-ver-sals treach-ery turn-around un-at-tached un-err-ing-ly
white-space wide-spread wing-spread wretch-ed wretch-ed-ly Brown-ian
Eng-lish Euler-ian Feb-ru-ary Gauss-ian Grothen-dieck Hamil-ton-ian
Her-mit-ian Jan-u-ary Japan-ese Kor-te-weg Le-gendre Lip-schitz
Lip-schitz-ian Mar-kov-ian Noe-ther-ian No-vem-ber Rie-mann-ian
Schwarz-schild Sep-tem-ber
%Zus"tze (Sind auch in besp.exc einzutragen!):
form per-iods Uni-ver-si-ty cri-ti-sism for-ma-lism}
\Invalid@\nofrills
\Invalid@\usualspace
\newif\ifnofrills@
\def\nofrills@#1#2{\relaxnext@
  \DN@{\ifx\next\nofrills
    \nofrills@true\let#2\relax\DN@\nofrills{\nextii@}%
  \else
    \nofrills@false\def#2{#1}\let\next@\nextii@\fi
\next@}}
\def\usualspace@#1{\ifnofrills@\def\usualspace{#1}\fi}
\def\addto#1#2{\csname \expandafter\eat@\string#1@\endcsname
  \expandafter{\the\csname \expandafter\eat@\string#1@\endcsname#2}}
\newdimen\bigsize@
\def\big@#1#2{{\hbox{$\left#2\vcenter to#1\bigsize@{}%
  \right.\nulldelimiterspace\z@\m@th$}}}
\def\big{\big@\@ne}
\def\Big{\big@{1.5}}
\def\bigg{\big@\tw@}
\def\Bigg{\big@{2.5}}
\def\raggedcenter@{\leftskip\z@ plus.4\hsize \rightskip\leftskip
 \parfillskip\z@ \parindent\z@ \spaceskip.3333em \xspaceskip.5em
 \pretolerance9999\tolerance9999 \exhyphenpenalty\@M
 \hyphenpenalty\@M \let\\\linebreak}
\def\upperspecialchars{\def\ss{SS}\let\i=I\let\j=J\let\ae\AE\let\oe\OE
  \let\o\O\let\aa\AA\let\l\L}
\def\uppercasetext@#1{%
  {\spaceskip1.2\fontdimen2\the\font plus1.2\fontdimen3\the\font
   \upperspecialchars\uctext@#1$\m@th\aftergroup\eat@$}}
\def\uctext@#1$#2${\endash@#1-\endash@$#2$\uctext@}
\def\endash@#1-#2\endash@{\uppercase{#1}\if\notempty{#2}--\endash@#2\endash@\fi}
\def\runaway@#1{\DN@{#1}\ifx\envir@\next@
  \Err@{You seem to have a missing or misspelled \string\end#1 ...}%
  \let\envir@\empty\fi}
\newif\iftemp@
\def\notempty#1{TT\fi\def\test@{#1}\ifx\test@\empty\temp@false
  \else\temp@true\fi \iftemp@}
\font@\tensmc=cmcsc10
\font@\sevenex=cmex7
\font@\sevenit=cmti7
\font@\eightrm=cmr8 % preloaded in plain.tex
\font@\sixrm=cmr6 % preloaded in plain.tex
\font@\eighti=cmmi8     \skewchar\eighti='177 % preloaded
\font@\sixi=cmmi6       \skewchar\sixi='177   % preloaded
\font@\eightsy=cmsy8    \skewchar\eightsy='60 % preloaded
\font@\sixsy=cmsy6      \skewchar\sixsy='60   % preloaded
\font@\eightex=cmex8
\font@\eightbf=cmbx8 % preloaded in plain.tex
\font@\sixbf=cmbx6   % preloaded in plain.tex
\font@\eightit=cmti8 % preloaded in plain.tex
\font@\eightsl=cmsl8 % preloaded in plain.tex
\font@\eightsmc=cmcsc8
\font@\eighttt=cmtt8 % preloaded in plain.tex
%\font@\ninerm=cmr9
%\font@\ninei=cmmi9    \skewchar\ninei='177
%\font@\ninesy=cmsy9   \skewchar\ninesy='60
%\font@\nineex=cmex9
%\font@\ninebf=cmbx9
%\font@\nineit=cmti9
%\font@\ninesl=cmsl9
%\font@\ninesmc=cmcsc9
%\font@\ninemsa=msam9
%\font@\ninemsb=msbm9
%\font@\nineeufm=eufm9

%Erg"nzung des fetten Small-Capitals-Fonts:
%\font@\eightbsmc=cmbcsc10 scaled 833
%\font@\tenbsmc=cmbcsc10
%\font@\twelvebsmc=cmbcsc10 scaled \magstep1
%\font@\fourteenbsmc=cmbcsc10 scaled \magstep2
%\font@\seventeenbsmc=cmbcsc10 scaled \magstep3
%\font@\twentybsmc=cmbcsc10 scaled \magstep4

\loadmsam
\loadmsbm
\loadeufm
\UseAMSsymbols
\newtoks\tenpoint@
\def\tenpoint{\normalbaselineskip12\p@
 \abovedisplayskip12\p@ plus3\p@ minus9\p@
 \belowdisplayskip\abovedisplayskip
 \abovedisplayshortskip\z@ plus3\p@
 \belowdisplayshortskip7\p@ plus3\p@ minus4\p@
 \textonlyfont@\rm\tenrm \textonlyfont@\it\tenit
 \textonlyfont@\sl\tensl \textonlyfont@\bf\tenbf
 \textonlyfont@\smc\tensmc \textonlyfont@\tt\tentt
%Erg"nzung des fetten Small-Capitals-Fonts:
 \textonlyfont@\bsmc\tenbsmc
 \ifsyntax@ \def\big##1{{\hbox{$\left##1\right.$}}}%
  \let\Big\big \let\bigg\big \let\Bigg\big
 \else
  \textfont\z@=\tenrm  \scriptfont\z@=\sevenrm  \scriptscriptfont\z@=\fiverm
  \textfont\@ne=\teni  \scriptfont\@ne=\seveni  \scriptscriptfont\@ne=\fivei
  \textfont\tw@=\tensy \scriptfont\tw@=\sevensy \scriptscriptfont\tw@=\fivesy
  \textfont\thr@@=\tenex \scriptfont\thr@@=\sevenex
        \scriptscriptfont\thr@@=\sevenex
  \textfont\itfam=\tenit \scriptfont\itfam=\sevenit
        \scriptscriptfont\itfam=\sevenit
  \textfont\bffam=\tenbf \scriptfont\bffam=\sevenbf
        \scriptscriptfont\bffam=\fivebf
  \setbox\strutbox\hbox{\vrule height8.5\p@ depth3.5\p@ width\z@}%
  \setbox\strutbox@\hbox{\lower.5\normallineskiplimit\vbox{%
        \kern-\normallineskiplimit\copy\strutbox}}%
 \setbox\z@\vbox{\hbox{$($}\kern\z@}\bigsize@=1.2\ht\z@
 \fi
 \normalbaselines\rm\ex@.2326ex\jot3\ex@\the\tenpoint@}
\newtoks\eightpoint@
\def\eightpoint{\normalbaselineskip10\p@
 \abovedisplayskip10\p@ plus2.4\p@ minus7.2\p@
 \belowdisplayskip\abovedisplayskip
 \abovedisplayshortskip\z@ plus2.4\p@
 \belowdisplayshortskip5.6\p@ plus2.4\p@ minus3.2\p@
 \textonlyfont@\rm\eightrm \textonlyfont@\it\eightit
 \textonlyfont@\sl\eightsl \textonlyfont@\bf\eightbf
 \textonlyfont@\smc\eightsmc \textonlyfont@\tt\eighttt
%Erg"nzung des fetten Small-Capitals-Fonts:
 \textonlyfont@\bsmc\eightbsmc
 \ifsyntax@\def\big##1{{\hbox{$\left##1\right.$}}}%
  \let\Big\big \let\bigg\big \let\Bigg\big
 \else
  \textfont\z@=\eightrm \scriptfont\z@=\sixrm \scriptscriptfont\z@=\fiverm
  \textfont\@ne=\eighti \scriptfont\@ne=\sixi \scriptscriptfont\@ne=\fivei
  \textfont\tw@=\eightsy \scriptfont\tw@=\sixsy \scriptscriptfont\tw@=\fivesy
  \textfont\thr@@=\eightex \scriptfont\thr@@=\sevenex
   \scriptscriptfont\thr@@=\sevenex
  \textfont\itfam=\eightit \scriptfont\itfam=\sevenit
   \scriptscriptfont\itfam=\sevenit
  \textfont\bffam=\eightbf \scriptfont\bffam=\sixbf
   \scriptscriptfont\bffam=\fivebf
 \setbox\strutbox\hbox{\vrule height7\p@ depth3\p@ width\z@}%
 \setbox\strutbox@\hbox{\raise.5\normallineskiplimit\vbox{%
   \kern-\normallineskiplimit\copy\strutbox}}%
 \setbox\z@\vbox{\hbox{$($}\kern\z@}\bigsize@=1.2\ht\z@
 \fi
 \normalbaselines\eightrm\ex@.2326ex\jot3\ex@\the\eightpoint@}

%Definition von 12-point, 14-point und 17-point Fonts
\font@\twelverm=cmr10 scaled\magstep1
\font@\twelveit=cmti10 scaled\magstep1
\font@\twelvesl=cmsl10 scaled\magstep1
\font@\twelvesmc=cmcsc10 scaled\magstep1
\font@\twelvett=cmtt10 scaled\magstep1
\font@\twelvebf=cmbx10 scaled\magstep1
\font@\twelvei=cmmi10 scaled\magstep1
\font@\twelvesy=cmsy10 scaled\magstep1
\font@\twelveex=cmex10 scaled\magstep1
\font@\twelvemsa=msam10 scaled\magstep1
\font@\twelveeufm=eufm10 scaled\magstep1
\font@\twelvemsb=msbm10 scaled\magstep1
\newtoks\twelvepoint@
\def\twelvepoint{\normalbaselineskip15\p@
 \abovedisplayskip15\p@ plus3.6\p@ minus10.8\p@
 \belowdisplayskip\abovedisplayskip
 \abovedisplayshortskip\z@ plus3.6\p@
 \belowdisplayshortskip8.4\p@ plus3.6\p@ minus4.8\p@
 \textonlyfont@\rm\twelverm \textonlyfont@\it\twelveit
 \textonlyfont@\sl\twelvesl \textonlyfont@\bf\twelvebf
 \textonlyfont@\smc\twelvesmc \textonlyfont@\tt\twelvett
%Erg"nzung des fetten Small-Capitals-Fonts:
 \textonlyfont@\bsmc\twelvebsmc
 \ifsyntax@ \def\big##1{{\hbox{$\left##1\right.$}}}%
  \let\Big\big \let\bigg\big \let\Bigg\big
 \else
  \textfont\z@=\twelverm  \scriptfont\z@=\tenrm  \scriptscriptfont\z@=\sevenrm
  \textfont\@ne=\twelvei  \scriptfont\@ne=\teni  \scriptscriptfont\@ne=\seveni
  \textfont\tw@=\twelvesy \scriptfont\tw@=\tensy \scriptscriptfont\tw@=\sevensy
  \textfont\thr@@=\twelveex \scriptfont\thr@@=\tenex
        \scriptscriptfont\thr@@=\tenex
  \textfont\itfam=\twelveit \scriptfont\itfam=\tenit
        \scriptscriptfont\itfam=\tenit
  \textfont\bffam=\twelvebf \scriptfont\bffam=\tenbf
        \scriptscriptfont\bffam=\sevenbf
  \setbox\strutbox\hbox{\vrule height10.2\p@ depth4.2\p@ width\z@}%
  \setbox\strutbox@\hbox{\lower.6\normallineskiplimit\vbox{%
        \kern-\normallineskiplimit\copy\strutbox}}%
 \setbox\z@\vbox{\hbox{$($}\kern\z@}\bigsize@=1.4\ht\z@
 \fi
 \normalbaselines\rm\ex@.2326ex\jot3.6\ex@\the\twelvepoint@}

\def\headfonts{\twelvepoint\bf}

\font@\fourteenrm=cmr10 scaled\magstep2
\font@\fourteenit=cmti10 scaled\magstep2
\font@\fourteensl=cmsl10 scaled\magstep2
\font@\fourteensmc=cmcsc10 scaled\magstep2
\font@\fourteentt=cmtt10 scaled\magstep2
\font@\fourteenbf=cmbx10 scaled\magstep2
\font@\fourteeni=cmmi10 scaled\magstep2
\font@\fourteensy=cmsy10 scaled\magstep2
\font@\fourteenex=cmex10 scaled\magstep2
\font@\fourteenmsa=msam10 scaled\magstep2
\font@\fourteeneufm=eufm10 scaled\magstep2
\font@\fourteenmsb=msbm10 scaled\magstep2
\newtoks\fourteenpoint@
\def\fourteenpoint{\normalbaselineskip15\p@
 \abovedisplayskip18\p@ plus4.3\p@ minus12.9\p@
 \belowdisplayskip\abovedisplayskip
 \abovedisplayshortskip\z@ plus4.3\p@
 \belowdisplayshortskip10.1\p@ plus4.3\p@ minus5.8\p@
 \textonlyfont@\rm\fourteenrm \textonlyfont@\it\fourteenit
 \textonlyfont@\sl\fourteensl \textonlyfont@\bf\fourteenbf
 \textonlyfont@\smc\fourteensmc \textonlyfont@\tt\fourteentt
%Erg"nzung des fetten Small-Capitals-Fonts:
 \textonlyfont@\bsmc\fourteenbsmc
 \ifsyntax@ \def\big##1{{\hbox{$\left##1\right.$}}}%
  \let\Big\big \let\bigg\big \let\Bigg\big
 \else
  \textfont\z@=\fourteenrm  \scriptfont\z@=\twelverm  \scriptscriptfont\z@=\tenrm
  \textfont\@ne=\fourteeni  \scriptfont\@ne=\twelvei  \scriptscriptfont\@ne=\teni
  \textfont\tw@=\fourteensy \scriptfont\tw@=\twelvesy \scriptscriptfont\tw@=\tensy
  \textfont\thr@@=\fourteenex \scriptfont\thr@@=\twelveex
        \scriptscriptfont\thr@@=\twelveex
  \textfont\itfam=\fourteenit \scriptfont\itfam=\twelveit
        \scriptscriptfont\itfam=\twelveit
  \textfont\bffam=\fourteenbf \scriptfont\bffam=\twelvebf
        \scriptscriptfont\bffam=\tenbf
  \setbox\strutbox\hbox{\vrule height12.2\p@ depth5\p@ width\z@}%
  \setbox\strutbox@\hbox{\lower.72\normallineskiplimit\vbox{%
        \kern-\normallineskiplimit\copy\strutbox}}%
 \setbox\z@\vbox{\hbox{$($}\kern\z@}\bigsize@=1.7\ht\z@
 \fi
 \normalbaselines\rm\ex@.2326ex\jot4.3\ex@\the\fourteenpoint@}

\def\chapheadfonts{\fourteenpoint\bf}

\font@\seventeenrm=cmr10 scaled\magstep3
\font@\seventeenit=cmti10 scaled\magstep3
\font@\seventeensl=cmsl10 scaled\magstep3
\font@\seventeensmc=cmcsc10 scaled\magstep3
\font@\seventeentt=cmtt10 scaled\magstep3
\font@\seventeenbf=cmbx10 scaled\magstep3
\font@\seventeeni=cmmi10 scaled\magstep3
\font@\seventeensy=cmsy10 scaled\magstep3
\font@\seventeenex=cmex10 scaled\magstep3
\font@\seventeenmsa=msam10 scaled\magstep3
\font@\seventeeneufm=eufm10 scaled\magstep3
\font@\seventeenmsb=msbm10 scaled\magstep3
\newtoks\seventeenpoint@
\def\seventeenpoint{\normalbaselineskip18\p@
 \abovedisplayskip21.6\p@ plus5.2\p@ minus15.4\p@
 \belowdisplayskip\abovedisplayskip
 \abovedisplayshortskip\z@ plus5.2\p@
 \belowdisplayshortskip12.1\p@ plus5.2\p@ minus7\p@
 \textonlyfont@\rm\seventeenrm \textonlyfont@\it\seventeenit
 \textonlyfont@\sl\seventeensl \textonlyfont@\bf\seventeenbf
 \textonlyfont@\smc\seventeensmc \textonlyfont@\tt\seventeentt
%Erg"nzung des fetten Small-Capitals-Fonts:
 \textonlyfont@\bsmc\seventeenbsmc
 \ifsyntax@ \def\big##1{{\hbox{$\left##1\right.$}}}%
  \let\Big\big \let\bigg\big \let\Bigg\big
 \else
  \textfont\z@=\seventeenrm  \scriptfont\z@=\fourteenrm  \scriptscriptfont\z@=\twelverm
  \textfont\@ne=\seventeeni  \scriptfont\@ne=\fourteeni  \scriptscriptfont\@ne=\twelvei
  \textfont\tw@=\seventeensy \scriptfont\tw@=\fourteensy \scriptscriptfont\tw@=\twelvesy
  \textfont\thr@@=\seventeenex \scriptfont\thr@@=\fourteenex
        \scriptscriptfont\thr@@=\fourteenex
  \textfont\itfam=\seventeenit \scriptfont\itfam=\fourteenit
        \scriptscriptfont\itfam=\fourteenit
  \textfont\bffam=\seventeenbf \scriptfont\bffam=\fourteenbf
        \scriptscriptfont\bffam=\twelvebf
  \setbox\strutbox\hbox{\vrule height14.6\p@ depth6\p@ width\z@}%
  \setbox\strutbox@\hbox{\lower.86\normallineskiplimit\vbox{%
        \kern-\normallineskiplimit\copy\strutbox}}%
 \setbox\z@\vbox{\hbox{$($}\kern\z@}\bigsize@=2\ht\z@
 \fi
 \normalbaselines\rm\ex@.2326ex\jot5.2\ex@\the\seventeenpoint@}

\font@\rrrrrm=cmr10 scaled\magstep4
\font@\bigtitlefont=cmbx10 scaled\magstep4

\parindent1pc
\normallineskiplimit\p@
\newdimen\indenti \indenti=2pc
\def\pageheight#1{\vsize#1}
\def\pagewidth#1{\hsize#1%
   \captionwidth@\hsize \advance\captionwidth@-2\indenti}
\pagewidth{30pc} \pageheight{47pc}
\def\topmatter{%
 \ifx\undefined\msafam
 \else\font@\eightmsa=msam8 \font@\sixmsa=msam6
   \ifsyntax@\else \addto\tenpoint{\textfont\msafam=\tenmsa
              \scriptfont\msafam=\sevenmsa \scriptscriptfont\msafam=\fivemsa}%
     \addto\eightpoint{\textfont\msafam=\eightmsa \scriptfont\msafam=\sixmsa
              \scriptscriptfont\msafam=\fivemsa}%
   \fi
 \fi
 \ifx\undefined\msbfam
 \else\font@\eightmsb=msbm8 \font@\sixmsb=msbm6
   \ifsyntax@\else \addto\tenpoint{\textfont\msbfam=\tenmsb
         \scriptfont\msbfam=\sevenmsb \scriptscriptfont\msbfam=\fivemsb}%
     \addto\eightpoint{\textfont\msbfam=\eightmsb \scriptfont\msbfam=\sixmsb
         \scriptscriptfont\msbfam=\fivemsb}%
   \fi
 \fi
 \ifx\undefined\eufmfam
 \else \font@\eighteufm=eufm8 \font@\sixeufm=eufm6
   \ifsyntax@\else \addto\tenpoint{\textfont\eufmfam=\teneufm
       \scriptfont\eufmfam=\seveneufm \scriptscriptfont\eufmfam=\fiveeufm}%
     \addto\eightpoint{\textfont\eufmfam=\eighteufm
       \scriptfont\eufmfam=\sixeufm \scriptscriptfont\eufmfam=\fiveeufm}%
   \fi
 \fi
 \ifx\undefined\eufbfam
 \else \font@\eighteufb=eufb8 \font@\sixeufb=eufb6
   \ifsyntax@\else \addto\tenpoint{\textfont\eufbfam=\teneufb
      \scriptfont\eufbfam=\seveneufb \scriptscriptfont\eufbfam=\fiveeufb}%
    \addto\eightpoint{\textfont\eufbfam=\eighteufb
      \scriptfont\eufbfam=\sixeufb \scriptscriptfont\eufbfam=\fiveeufb}%
   \fi
 \fi
 \ifx\undefined\eusmfam
 \else \font@\eighteusm=eusm8 \font@\sixeusm=eusm6
   \ifsyntax@\else \addto\tenpoint{\textfont\eusmfam=\teneusm
       \scriptfont\eusmfam=\seveneusm \scriptscriptfont\eusmfam=\fiveeusm}%
     \addto\eightpoint{\textfont\eusmfam=\eighteusm
       \scriptfont\eusmfam=\sixeusm \scriptscriptfont\eusmfam=\fiveeusm}%
   \fi
 \fi
 \ifx\undefined\eusbfam
 \else \font@\eighteusb=eusb8 \font@\sixeusb=eusb6
   \ifsyntax@\else \addto\tenpoint{\textfont\eusbfam=\teneusb
       \scriptfont\eusbfam=\seveneusb \scriptscriptfont\eusbfam=\fiveeusb}%
     \addto\eightpoint{\textfont\eusbfam=\eighteusb
       \scriptfont\eusbfam=\sixeusb \scriptscriptfont\eusbfam=\fiveeusb}%
   \fi
 \fi
 \ifx\undefined\eurmfam
 \else \font@\eighteurm=eurm8 \font@\sixeurm=eurm6
   \ifsyntax@\else \addto\tenpoint{\textfont\eurmfam=\teneurm
       \scriptfont\eurmfam=\seveneurm \scriptscriptfont\eurmfam=\fiveeurm}%
     \addto\eightpoint{\textfont\eurmfam=\eighteurm
       \scriptfont\eurmfam=\sixeurm \scriptscriptfont\eurmfam=\fiveeurm}%
   \fi
 \fi
 \ifx\undefined\eurbfam
 \else \font@\eighteurb=eurb8 \font@\sixeurb=eurb6
   \ifsyntax@\else \addto\tenpoint{\textfont\eurbfam=\teneurb
       \scriptfont\eurbfam=\seveneurb \scriptscriptfont\eurbfam=\fiveeurb}%
    \addto\eightpoint{\textfont\eurbfam=\eighteurb
       \scriptfont\eurbfam=\sixeurb \scriptscriptfont\eurbfam=\fiveeurb}%
   \fi
 \fi
 \ifx\undefined\cmmibfam
 \else \font@\eightcmmib=cmmib8 \font@\sixcmmib=cmmib6
   \ifsyntax@\else \addto\tenpoint{\textfont\cmmibfam=\tencmmib
       \scriptfont\cmmibfam=\sevencmmib \scriptscriptfont\cmmibfam=\fivecmmib}%
    \addto\eightpoint{\textfont\cmmibfam=\eightcmmib
       \scriptfont\cmmibfam=\sixcmmib \scriptscriptfont\cmmibfam=\fivecmmib}%
   \fi
 \fi
 \ifx\undefined\cmbsyfam
 \else \font@\eightcmbsy=cmbsy8 \font@\sixcmbsy=cmbsy6
   \ifsyntax@\else \addto\tenpoint{\textfont\cmbsyfam=\tencmbsy
      \scriptfont\cmbsyfam=\sevencmbsy \scriptscriptfont\cmbsyfam=\fivecmbsy}%
    \addto\eightpoint{\textfont\cmbsyfam=\eightcmbsy
      \scriptfont\cmbsyfam=\sixcmbsy \scriptscriptfont\cmbsyfam=\fivecmbsy}%
   \fi
 \fi
 \let\topmatter\relax}
\def\chapterno@{\uppercase\expandafter{\romannumeral\chaptercount@}}
\newcount\chaptercount@
\def\chapter{\nofrills@{\afterassignment\chapterno@
                        CHAPTER \global\chaptercount@=}\chapter@
 \DNii@##1{\leavevmode\hskip-\leftskip
   \rlap{\vbox to\z@{\vss\centerline{\eightpoint
   \chapter@##1\unskip}\baselineskip2pc\null}}\hskip\leftskip
   \nofrills@false}%
 \FN@\next@}
\newbox\titlebox@

%Uppercase ist abgestellt.
\def\title{\nofrills@{\relax}\title@%
 \DNii@##1\endtitle{\global\setbox\titlebox@\vtop{\tenpoint\bf
 \raggedcenter@\ignorespaces
 \baselineskip1.3\baselineskip\title@{##1}\endgraf}%
 \ifmonograph@ \edef\next{\the\leftheadtoks}\ifx\next\empty
    \leftheadtext{##1}\fi
 \fi
 \edef\next{\the\rightheadtoks}\ifx\next\empty \rightheadtext{##1}\fi
 }\FN@\next@}
\newbox\authorbox@
\def\author#1\endauthor{\global\setbox\authorbox@
 \vbox{\tenpoint\smc\raggedcenter@\ignorespaces
 #1\endgraf}\relaxnext@ \edef\next{\the\leftheadtoks}%
 \ifx\next\empty\leftheadtext{#1}\fi}
\newbox\affilbox@
\def\affil#1\endaffil{\global\setbox\affilbox@
 \vbox{\tenpoint\raggedcenter@\ignorespaces#1\endgraf}}
\newcount\addresscount@
\addresscount@\z@
\def\address#1\endaddress{\global\advance\addresscount@\@ne
  \expandafter\gdef\csname address\number\addresscount@\endcsname
  {\vskip12\p@ minus6\p@\noindent\eightpoint\smc\ignorespaces#1\par}}
\def\email{\nofrills@{\eightpoint{\it E-mail\/}:\enspace}\email@
  \DNii@##1\endemail{%
  \expandafter\gdef\csname email\number\addresscount@\endcsname
  {\def\usualspace{{\it\enspace}}\smallskip\noindent\eightpoint\email@
  \ignorespaces##1\par}}%
 \FN@\next@}
\def\thedate@{}
\def\date#1\enddate{\gdef\thedate@{\tenpoint\ignorespaces#1\unskip}}
\def\thethanks@{}
\def\thanks#1\endthanks{\gdef\thethanks@{\eightpoint\ignorespaces#1.\unskip}}
\def\thekeywords@{}
\def\keywords{\nofrills@{{\it Key words and phrases.\enspace}}\keywords@
 \DNii@##1\endkeywords{\def\thekeywords@{\def\usualspace{{\it\enspace}}%
 \eightpoint\keywords@\ignorespaces##1\unskip.}}%
 \FN@\next@}
\def\thesubjclass@{}
\def\subjclass{\nofrills@{{\rm2000 {\it Mathematics Subject
   Classification\/}.\enspace}}\subjclass@
 \DNii@##1\endsubjclass{\def\thesubjclass@{\def\usualspace
  {{\rm\enspace}}\eightpoint\subjclass@\ignorespaces##1\unskip.}}%
 \FN@\next@}
\newbox\abstractbox@
\def\abstract{\nofrills@{{\smc Abstract.\enspace}}\abstract@
 \DNii@{\setbox\abstractbox@\vbox\bgroup\noindent$$\vbox\bgroup
  \def\envir@{abstract}\advance\hsize-2\indenti
  \usualspace@{{\enspace}}\eightpoint \noindent\abstract@\ignorespaces}%
 \FN@\next@}
\def\endabstract{\par\unskip\egroup$$\egroup}
\def\widestnumber#1#2{\begingroup\let\head\null\let\subhead\empty
   \let\subsubhead\subhead
   \ifx#1\head\global\setbox\tocheadbox@\hbox{#2.\enspace}%
   \else\ifx#1\subhead\global\setbox\tocsubheadbox@\hbox{#2.\enspace}%
   \else\ifx#1\key\bgroup\let\endrefitem@\egroup
        \key#2\endrefitem@\global\refindentwd\wd\keybox@
   \else\ifx#1\no\bgroup\let\endrefitem@\egroup
        \no#2\endrefitem@\global\refindentwd\wd\nobox@
   \else\ifx#1\page\global\setbox\pagesbox@\hbox{\quad\bf#2}%
   \else\ifx#1\item\setboxz@h{#2}\global\rosteritemwd\wdz@
        \global\advance\rosteritemwd by.5\parindent
   \else\message{\string\widestnumber is not defined for this option
   (\string#1)}%
\fi\fi\fi\fi\fi\fi\endgroup}
\newif\ifmonograph@
\def\Monograph{\monograph@true \let\headmark\rightheadtext
  \let\varindent@\indent \def\headfont@{\bf}\def\proclaimheadfont@{\smc}%
  \def\demofont@{\smc}}
%Bei Proclaim,...: Einrcken.
\let\varindent@\indent

\newbox\tocheadbox@    \newbox\tocsubheadbox@
\newbox\tocbox@
\def\toc{\toc@{Contents}}
\def\newtocdefs{%
   \def \title##1\endtitle
       {\penaltyandskip@\z@\smallskipamount
        \hangindent\wd\tocheadbox@\noindent{\bf##1}}%
   \def \chapter##1{%
        Chapter \uppercase\expandafter{\romannumeral##1.\unskip}\enspace}%
   \def \specialhead##1\endspecialhead
       {\par\hangindent\wd\tocheadbox@ \noindent##1\par}%
   \def \head##1 ##2\endhead
       {\par\hangindent\wd\tocheadbox@ \noindent
        \if\notempty{##1}\hbox to\wd\tocheadbox@{\hfil##1\enspace}\fi
        ##2\par}%
   \def \subhead##1 ##2\endsubhead
       {\par\vskip-\parskip {\normalbaselines
        \advance\leftskip\wd\tocheadbox@
        \hangindent\wd\tocsubheadbox@ \noindent
        \if\notempty{##1}\hbox to\wd\tocsubheadbox@{##1\unskip\hfil}\fi
         ##2\par}}%
   \def \subsubhead##1 ##2\endsubsubhead
       {\par\vskip-\parskip {\normalbaselines
        \advance\leftskip\wd\tocheadbox@
        \hangindent\wd\tocsubheadbox@ \noindent
        \if\notempty{##1}\hbox to\wd\tocsubheadbox@{##1\unskip\hfil}\fi
        ##2\par}}}
\def\toc@#1{\relaxnext@
   \def\page##1%
       {\unskip\penalty0\null\hfil
        \rlap{\hbox to\wd\pagesbox@{\quad\hfil##1}}\hfilneg\penalty\@M}%
 \DN@{\ifx\next\nofrills\DN@\nofrills{\nextii@}%
      \else\DN@{\nextii@{{#1}}}\fi
      \next@}%
 \DNii@##1{%
\ifmonograph@\bgroup\else\setbox\tocbox@\vbox\bgroup
   \centerline{\headfont@\ignorespaces##1\unskip}\nobreak
   \vskip\belowheadskip \fi
   \setbox\tocheadbox@\hbox{0.\enspace}%
   \setbox\tocsubheadbox@\hbox{0.0.\enspace}%
   \leftskip\indenti \rightskip\leftskip
   \setbox\pagesbox@\hbox{\bf\quad000}\advance\rightskip\wd\pagesbox@
   \newtocdefs
 }%
 \FN@\next@}
\def\endtoc{\par\egroup}
\let\pretitle\relax
\let\preauthor\relax
\let\preaffil\relax
\let\predate\relax
\let\preabstract\relax
\let\prepaper\relax
\def\dedicatory #1\enddedicatory{\def\preabstract{{\medskip
  \eightpoint\it \raggedcenter@#1\endgraf}}}
\def\thetranslator@{}
\def\translator#1\endtranslator{\def\thetranslator@{\nobreak\medskip
 \line{\eightpoint\hfil Translated by \uppercase{#1}\qquad\qquad}\nobreak}}
\outer\def\endtopmatter{\runaway@{abstract}%
 \edef\next{\the\leftheadtoks}\ifx\next\empty
  \expandafter\leftheadtext\expandafter{\the\rightheadtoks}\fi
 \ifmonograph@\else
   \ifx\thesubjclass@\empty\else \makefootnote@{}{\thesubjclass@}\fi
   \ifx\thekeywords@\empty\else \makefootnote@{}{\thekeywords@}\fi
   \ifx\thethanks@\empty\else \makefootnote@{}{\thethanks@}\fi
 \fi
  \pretitle
  \ifmonograph@ \topskip7pc \else \topskip4pc \fi
  \box\titlebox@
  \topskip10pt% reset to normal value
  \preauthor
  \ifvoid\authorbox@\else \vskip2.5pc plus1pc \unvbox\authorbox@\fi
  \preaffil
  \ifvoid\affilbox@\else \vskip1pc plus.5pc \unvbox\affilbox@\fi
  \predate
  \ifx\thedate@\empty\else \vskip1pc plus.5pc \line{\hfil\thedate@\hfil}\fi
  \preabstract
  \ifvoid\abstractbox@\else \vskip1.5pc plus.5pc \unvbox\abstractbox@ \fi
  \ifvoid\tocbox@\else\vskip1.5pc plus.5pc \unvbox\tocbox@\fi
  \prepaper
  \vskip2pc plus1pc
}
\def\document{\let\fontlist@\relax\let\alloclist@\relax
  \tenpoint}

%Modifizierte Head-Skips
\newskip\aboveheadskip       \aboveheadskip1.8\bigskipamount
\newdimen\belowheadskip      \belowheadskip1.8\medskipamount

\def\headfont@{\smc}
\def\penaltyandskip@#1#2{\relax\ifdim\lastskip<#2\relax\removelastskip
      \ifnum#1=\z@\else\penalty@#1\relax\fi\vskip#2%
  \else\ifnum#1=\z@\else\penalty@#1\relax\fi\fi}
\def\nobreak{\penalty\@M
  \ifvmode\def\penalty@{\let\penalty@\penalty\count@@@}%
  \everypar{\let\penalty@\penalty\everypar{}}\fi}
\let\penalty@\penalty
\def\heading#1\endheading{\head#1\endhead}

\def\specialheadfont@{\bf}
\outer\def\specialhead{\par\penaltyandskip@{-200}\aboveheadskip
  \begingroup\interlinepenalty\@M\rightskip\z@ plus\hsize \let\\\linebreak
  \specialheadfont@\noindent\ignorespaces}
\def\endspecialhead{\par\endgroup\nobreak\vskip\belowheadskip}
%\outer\def\head#1\endhead{\par\penaltyandskip@{-200}\aboveheadskip
% {\headfont@\raggedcenter@\interlinepenalty\@M
% \ignorespaces#1\endgraf}\nobreak
% \vskip\belowheadskip
% \headmark{#1}}
\let\headmark\eat@
\newskip\subheadskip       \subheadskip\medskipamount
\def\subheadfont@{\bf}
\outer\def\subhead{\nofrills@{.\enspace}\subhead@
 \DNii@##1\endsubhead{\par\penaltyandskip@{-100}\subheadskip
  \varindent@{\usualspace@{{\subheadfont@\enspace}}%
 \subheadfont@\ignorespaces##1\unskip\subhead@}\ignorespaces}%
 \FN@\next@}
\outer\def\subsubhead{\nofrills@{.\enspace}\subsubhead@
 \DNii@##1\endsubsubhead{\par\penaltyandskip@{-50}\medskipamount
      {\usualspace@{{\it\enspace}}%
  \it\ignorespaces##1\unskip\subsubhead@}\ignorespaces}%
 \FN@\next@}
\def\proclaimheadfont@{\bf}
\outer\def\proclaim{\runaway@{proclaim}\def\envir@{proclaim}%
  \nofrills@{.\enspace}\proclaim@
 \DNii@##1{\penaltyandskip@{-100}\medskipamount\varindent@
   \usualspace@{{\proclaimheadfont@\enspace}}\proclaimheadfont@
   \ignorespaces##1\unskip\proclaim@
  \sl\ignorespaces}% 
 \FN@\next@}
\outer\def\endproclaim{\let\envir@\relax\par\rm
  \penaltyandskip@{55}\medskipamount}
\def\demoheadfont@{\it}
\def\demo{\runaway@{proclaim}\nofrills@{.\enspace}\demo@
     \DNii@##1{\par\penaltyandskip@\z@\medskipamount
  {\usualspace@{{\demoheadfont@\enspace}}%
  \varindent@\demoheadfont@\ignorespaces##1\unskip\demo@}\rm
  \ignorespaces}\FN@\next@}
\def\enddemo{\par\medskip}
\def\qed{\ifhmode\unskip\nobreak\fi\quad\ifmmode\square\else$\m@th\square$\fi}

%\endremark=\enddemo

%Headfont fr Definition wie bei Beweisen.
\def\definition{\runaway@{proclaim}%
  \nofrills@{.\demoheadfont@\enspace}\definition@
        \DNii@##1{\penaltyandskip@{-100}\medskipamount
        {\usualspace@{{\demoheadfont@\enspace}}%
        \varindent@\demoheadfont@\ignorespaces##1\unskip\definition@}%
        \rm \ignorespaces}\FN@\next@}

%\let\example\definition
%\let\endexample\enddefinition
%Modifikation:

\newdimen\rosteritemwd
\newcount\rostercount@
\newif\iffirstitem@
\let\plainitem@\item
\newtoks\everypartoks@
\def\par@{\everypartoks@\expandafter{\the\everypar}\everypar{}}
\def\roster{\edef\leftskip@{\leftskip\the\leftskip}%
 \relaxnext@
 \rostercount@\z@  
 \def\item{\FN@\rosteritem@}% 
 \DN@{\ifx\next\runinitem\let\next@\nextii@\else
  \let\next@\nextiii@\fi\next@}%
 \DNii@\runinitem  
  {\unskip  
   \DN@{\ifx\next[\let\next@\nextii@\else
    \ifx\next"\let\next@\nextiii@\else\let\next@\nextiv@\fi\fi\next@}%
   \DNii@[####1]{\rostercount@####1\relax
    \enspace{\rm(\number\rostercount@)}~\ignorespaces}%
   \def\nextiii@"####1"{\enspace{\rm####1}~\ignorespaces}%
   \def\nextiv@{\enspace{\rm(1)}\rostercount@\@ne~}%
   \par@\firstitem@false  
   \FN@\next@}% 
 \def\nextiii@{\par\par@  
  \penalty\@m\smallskip\vskip-\parskip
  \firstitem@true}%  
 \FN@\next@}
\def\rosteritem@{\iffirstitem@\firstitem@false\else\par\vskip-\parskip\fi
 \leftskip3\parindent\noindent  
 \DNii@[##1]{\rostercount@##1\relax
  \llap{\hbox to2.5\parindent{\hss\rm(\number\rostercount@)}%
   \hskip.5\parindent}\ignorespaces}%
 \def\nextiii@"##1"{%
  \llap{\hbox to2.5\parindent{\hss\rm##1}\hskip.5\parindent}\ignorespaces}%
 \def\nextiv@{\advance\rostercount@\@ne
  \llap{\hbox to2.5\parindent{\hss\rm(\number\rostercount@)}%
   \hskip.5\parindent}}%
 \ifx\next[\let\next@\nextii@\else\ifx\next"\let\next@\nextiii@\else
  \let\next@\nextiv@\fi\fi\next@}

\newif\ifnextRunin@
\def\endroster{\relaxnext@
 \par\leftskip@  
 \penalty-50 \vskip-\parskip\smallskip  
 \DN@{\ifx\next\Runinitem\let\next@\relax
  \else\nextRunin@false\let\item\plainitem@  
   \ifx\next\par 
    \DN@\par{\everypar\expandafter{\the\everypartoks@}}%
   \else  
    \DN@{\noindent\everypar\expandafter{\the\everypartoks@}}%
  \fi\fi\next@}%
 \FN@\next@}
\newcount\rosterhangafter@
\def\Runinitem#1\roster\runinitem{\relaxnext@
 \rostercount@\z@ 
 \def\item{\FN@\rosteritem@}%  
 \def\runinitem@{#1}%  
 \DN@{\ifx\next[\let\next\nextii@\else\ifx\next"\let\next\nextiii@
  \else\let\next\nextiv@\fi\fi\next}%
 \DNii@[##1]{\rostercount@##1\relax
  \def\item@{{\rm(\number\rostercount@)}}\nextv@}%
 \def\nextiii@"##1"{\def\item@{{\rm##1}}\nextv@}%
 \def\nextiv@{\advance\rostercount@\@ne
  \def\item@{{\rm(\number\rostercount@)}}\nextv@}%
 \def\nextv@{\setbox\z@\vbox  
  {\ifnextRunin@\noindent\fi  
  \runinitem@\unskip\enspace\item@~\par  
  \global\rosterhangafter@\prevgraf}% 
  \firstitem@false  
  \ifnextRunin@\else\par\fi  
  \hangafter\rosterhangafter@\hangindent3\parindent
  \ifnextRunin@\noindent\fi  
  \runinitem@\unskip\enspace 
  \item@~\ifnextRunin@\else\par@\fi  
  \nextRunin@true\ignorespaces}%  
 \FN@\next@}
\def\footmarkform@#1{$\m@th^{#1}$}
\let\thefootnotemark\footmarkform@
\def\makefootnote@#1#2{\insert\footins
 {\interlinepenalty\interfootnotelinepenalty
 \eightpoint\splittopskip\ht\strutbox\splitmaxdepth\dp\strutbox
 \floatingpenalty\@MM\leftskip\z@\rightskip\z@\spaceskip\z@\xspaceskip\z@
 \leavevmode{#1}\footstrut\ignorespaces#2\unskip\lower\dp\strutbox
 \vbox to\dp\strutbox{}}}
\newcount\footmarkcount@
\footmarkcount@\z@
\def\footnotemark{\let\@sf\empty\relaxnext@
 \ifhmode\edef\@sf{\spacefactor\the\spacefactor}\/\fi
 \DN@{\ifx[\next\let\next@\nextii@\else
  \ifx"\next\let\next@\nextiii@\else
  \let\next@\nextiv@\fi\fi\next@}%
 \DNii@[##1]{\footmarkform@{##1}\@sf}%
 \def\nextiii@"##1"{{##1}\@sf}%
 \def\nextiv@{\iffirstchoice@\global\advance\footmarkcount@\@ne\fi
  \footmarkform@{\number\footmarkcount@}\@sf}%
 \FN@\next@}
\def\footnotetext{\relaxnext@
 \DN@{\ifx[\next\let\next@\nextii@\else
  \ifx"\next\let\next@\nextiii@\else
  \let\next@\nextiv@\fi\fi\next@}%
 \DNii@[##1]##2{\makefootnote@{\footmarkform@{##1}}{##2}}%
 \def\nextiii@"##1"##2{\makefootnote@{##1}{##2}}%
 \def\nextiv@##1{\makefootnote@{\footmarkform@{\number\footmarkcount@}}{##1}}%
 \FN@\next@}
\def\footnote{\let\@sf\empty\relaxnext@
 \ifhmode\edef\@sf{\spacefactor\the\spacefactor}\/\fi
 \DN@{\ifx[\next\let\next@\nextii@\else
  \ifx"\next\let\next@\nextiii@\else
  \let\next@\nextiv@\fi\fi\next@}%
 \DNii@[##1]##2{\footnotemark[##1]\footnotetext[##1]{##2}}%
 \def\nextiii@"##1"##2{\footnotemark"##1"\footnotetext"##1"{##2}}%
 \def\nextiv@##1{\footnotemark\footnotetext{##1}}%
 \FN@\next@}
\def\adjustfootnotemark#1{\advance\footmarkcount@#1\relax}
\def\footnoterule{\kern-3\p@
  \hrule width 5pc\kern 2.6\p@} 
\def\captionfont@{\smc}
\def\topcaption#1#2\endcaption{%
  {\dimen@\hsize \advance\dimen@-\captionwidth@
   \rm\raggedcenter@ \advance\leftskip.5\dimen@ \rightskip\leftskip
  {\captionfont@#1}%
  \if\notempty{#2}.\enspace\ignorespaces#2\fi
  \endgraf}\nobreak\bigskip}
\def\botcaption#1#2\endcaption{%
  \nobreak\bigskip
  \setboxz@h{\captionfont@#1\if\notempty{#2}.\enspace\rm#2\fi}%
  {\dimen@\hsize \advance\dimen@-\captionwidth@
   \leftskip.5\dimen@ \rightskip\leftskip
   \noindent \ifdim\wdz@>\captionwidth@ 
   \else\hfil\fi 
  {\captionfont@#1}\if\notempty{#2}.\enspace\rm#2\fi\endgraf}}
\def\@ins{\par\begingroup\def\vspace##1{\vskip##1\relax}%
  \def\captionwidth##1{\captionwidth@##1\relax}%
  \setbox\z@\vbox\bgroup} % start a \vbox
\def\block{\RIfMIfI@\nondmatherr@\block\fi
       \else\ifvmode\vskip\abovedisplayskip\noindent\fi
        $$\def\endblock{\par\egroup$$}\fi
  \vbox\bgroup\advance\hsize-2\indenti\noindent}
\def\endblock{\par\egroup}
\def\cite#1{{\rm[{\citefont@\m@th#1}]}}
\def\citefont@{\rm}
\def\refsfont@{\eightpoint}
\outer\def\Refs{\runaway@{proclaim}%
 \relaxnext@ \DN@{\ifx\next\nofrills\DN@\nofrills{\nextii@}\else
  \DN@{\nextii@{References}}\fi\next@}%
 \DNii@##1{\penaltyandskip@{-200}\aboveheadskip
  \line{\hfil\headfont@\ignorespaces##1\unskip\hfil}\nobreak
  \vskip\belowheadskip
  \begingroup\refsfont@\sfcode`.=\@m}%
 \FN@\next@}
\def\endRefs{\par\endgroup}
\newbox\nobox@            \newbox\keybox@           \newbox\bybox@
\newbox\paperbox@         \newbox\paperinfobox@     \newbox\jourbox@
\newbox\volbox@           \newbox\issuebox@         \newbox\yrbox@
\newbox\pagesbox@         \newbox\bookbox@          \newbox\bookinfobox@
\newbox\publbox@          \newbox\publaddrbox@      \newbox\finalinfobox@
\newbox\edsbox@           \newbox\langbox@
\newif\iffirstref@        \newif\iflastref@
\newif\ifprevjour@        \newif\ifbook@            \newif\ifprevinbook@
\newif\ifquotes@          \newif\ifbookquotes@      \newif\ifpaperquotes@
\newdimen\bysamerulewd@
\setboxz@h{\refsfont@\kern3em}
\bysamerulewd@\wdz@
\newdimen\refindentwd
\setboxz@h{\refsfont@ 00. }
\refindentwd\wdz@
\outer\def\ref{\begingroup \noindent\hangindent\refindentwd
 \firstref@true \def\nofrills{\def\refkern@{\kern3sp}}%
 \ref@}
\def\ref@{\book@false \bgroup\let\endrefitem@\egroup \ignorespaces}
\def\moreref{\endrefitem@\endref@\firstref@false\ref@}%
\def\transl{\endrefitem@\endref@\firstref@false
  \book@false
  \prepunct@
  \setboxz@h\bgroup \aftergroup\unhbox\aftergroup\z@
    \def\endrefitem@{\unskip\refkern@\egroup}\ignorespaces}%
\def\emptyifempty@{\dimen@\wd\currbox@
  \advance\dimen@-\wd\z@ \advance\dimen@-.1\p@
  \ifdim\dimen@<\z@ \setbox\currbox@\copy\voidb@x \fi}
\let\refkern@\relax
\def\endrefitem@{\unskip\refkern@\egroup
  \setboxz@h{\refkern@}\emptyifempty@}\ignorespaces
\def\refdef@#1#2#3{\edef\next@{\noexpand\endrefitem@
  \let\noexpand\currbox@\csname\expandafter\eat@\string#1box@\endcsname
    \noexpand\setbox\noexpand\currbox@\hbox\bgroup}%
  \toks@\expandafter{\next@}%
  \if\notempty{#2#3}\toks@\expandafter{\the\toks@
  \def\endrefitem@{\unskip#3\refkern@\egroup
  \setboxz@h{#2#3\refkern@}\emptyifempty@}#2}\fi
  \toks@\expandafter{\the\toks@\ignorespaces}%
  \edef#1{\the\toks@}}
\refdef@\no{}{. }
\refdef@\key{[\m@th}{] }
\refdef@\by{}{}
\def\bysame{\by\hbox to\bysamerulewd@{\hrulefill}\thinspace
   \kern0sp}
\def\manyby{\message{\string\manyby is no longer necessary; \string\by
  can be used instead, starting with version 2.0 of \styname.STY}\by}
\refdef@\paper{\ifpaperquotes@``\fi\it}{}
\refdef@\paperinfo{}{}
\def\jour{\endrefitem@\let\currbox@\jourbox@
  \setbox\currbox@\hbox\bgroup
  \def\endrefitem@{\unskip\refkern@\egroup
    \setboxz@h{\refkern@}\emptyifempty@
    \ifvoid\jourbox@\else\prevjour@true\fi}%
\ignorespaces}
\refdef@\vol{\ifbook@\else\bf\fi}{}
\refdef@\issue{no. }{}
\refdef@\yr{}{}
\refdef@\pages{}{}
\def\page{\endrefitem@\def\pp@{\def\pp@{pp.~}p.~}\let\currbox@\pagesbox@
  \setbox\currbox@\hbox\bgroup\ignorespaces}
\def\pp@{pp.~}
\def\book{\endrefitem@ \let\currbox@\bookbox@
 \setbox\currbox@\hbox\bgroup\def\endrefitem@{\unskip\refkern@\egroup
  \setboxz@h{\ifbookquotes@``\fi}\emptyifempty@
  \ifvoid\bookbox@\else\book@true\fi}%
  \ifbookquotes@``\fi\it\ignorespaces}
\def\inbook{\endrefitem@
  \let\currbox@\bookbox@\setbox\currbox@\hbox\bgroup
  \def\endrefitem@{\unskip\refkern@\egroup
  \setboxz@h{\ifbookquotes@``\fi}\emptyifempty@
  \ifvoid\bookbox@\else\book@true\previnbook@true\fi}%
  \ifbookquotes@``\fi\ignorespaces}
\refdef@\eds{(}{, eds.)}
\def\ed{\endrefitem@\let\currbox@\edsbox@
 \setbox\currbox@\hbox\bgroup
 \def\endrefitem@{\unskip, ed.)\refkern@\egroup
  \setboxz@h{(, ed.)}\emptyifempty@}(\ignorespaces}
\refdef@\bookinfo{}{}
\refdef@\publ{}{}
\refdef@\publaddr{}{}
\refdef@\finalinfo{}{}
\refdef@\lang{(}{)}

\let\refdef@\relax 
\def\ppunbox@#1{\ifvoid#1\else\prepunct@\unhbox#1\fi}
\def\nocomma@#1{\ifvoid#1\else\changepunct@3\prepunct@\unhbox#1\fi}
\def\changepunct@#1{\ifnum\lastkern<3 \unkern\kern#1sp\fi}
\def\prepunct@{\count@\lastkern\unkern
  \ifnum\lastpenalty=0
    \let\penalty@\relax
  \else
    \edef\penalty@{\penalty\the\lastpenalty\relax}%
  \fi
  \unpenalty
  \let\refspace@\ \ifcase\count@,% usual case, do a comma
\or;\or.\or % do nothing; this case is from nofrills.
  \or\let\refspace@\relax
  \else,\fi
  \ifquotes@''\quotes@false\fi \penalty@ \refspace@
}
\def\transferpenalty@#1{\dimen@\lastkern\unkern
  \ifnum\lastpenalty=0\unpenalty\let\penalty@\relax
  \else\edef\penalty@{\penalty\the\lastpenalty\relax}\unpenalty\fi
  #1\penalty@\kern\dimen@}
\def\endref{\endrefitem@\lastref@true\endref@
  \par\endgroup \prevjour@false \previnbook@false }
\def\endref@{%
\iffirstref@
  \ifvoid\nobox@\ifvoid\keybox@\indent\fi
  \else\hbox to\refindentwd{\hss\unhbox\nobox@}\fi
  \ifvoid\keybox@
  \else\ifdim\wd\keybox@>\refindentwd
         \box\keybox@
       \else\hbox to\refindentwd{\unhbox\keybox@\hfil}\fi\fi
  \kern4sp\ppunbox@\bybox@
\fi 
  \ifvoid\paperbox@
  \else\prepunct@\unhbox\paperbox@
    \ifpaperquotes@\quotes@true\fi\fi
  \ppunbox@\paperinfobox@
  \ifvoid\jourbox@
    \ifprevjour@ \nocomma@\volbox@
      \nocomma@\issuebox@
      \ifvoid\yrbox@\else\changepunct@3\prepunct@(\unhbox\yrbox@
        \transferpenalty@)\fi
      \ppunbox@\pagesbox@
    \fi 
  \else \prepunct@\unhbox\jourbox@
    \nocomma@\volbox@
    \nocomma@\issuebox@
    \ifvoid\yrbox@\else\changepunct@3\prepunct@(\unhbox\yrbox@
      \transferpenalty@)\fi
    \ppunbox@\pagesbox@
  \fi 
  \ifbook@\prepunct@\unhbox\bookbox@ \ifbookquotes@\quotes@true\fi \fi
  \nocomma@\edsbox@
  \ppunbox@\bookinfobox@
  \ifbook@\ifvoid\volbox@\else\prepunct@ vol.~\unhbox\volbox@
  \fi\fi
  \ppunbox@\publbox@ \ppunbox@\publaddrbox@
  \ifbook@ \ppunbox@\yrbox@
    \ifvoid\pagesbox@
    \else\prepunct@\pp@\unhbox\pagesbox@\fi
  \else
    \ifprevinbook@ \ppunbox@\yrbox@
      \ifvoid\pagesbox@\else\prepunct@\pp@\unhbox\pagesbox@\fi
    \fi \fi
  \ppunbox@\finalinfobox@
  \iflastref@
    \ifvoid\langbox@.\ifquotes@''\fi
    \else\changepunct@2\prepunct@\unhbox\langbox@\fi
  \else
    \ifvoid\langbox@\changepunct@1%
    \else\changepunct@3\prepunct@\unhbox\langbox@
      \changepunct@1\fi
  \fi
}
\outer\def\enddocument{%
 \runaway@{proclaim}%
\ifmonograph@ % do nothing
\else
 \nobreak
 \thetranslator@
 \count@\z@ \loop\ifnum\count@<\addresscount@\advance\count@\@ne
 \csname address\number\count@\endcsname
 \csname email\number\count@\endcsname
 \repeat
\fi
 \vfill\supereject\end}

%Modifizierte Fonts fr Proclaim,...
\def\headfont@{\headfonts}
\def\proclaimheadfont@{\bf}
\def\specialheadfont@{\bf}
\def\subheadfont@{\bf}
\def\demoheadfont@{\smc}

%Kontrollsequenzen fr Inhaltsverzeichnis und Index:
\newif\ifThisToToc \ThisToTocfalse
\newif\iftocloaded \tocloadedfalse

\def\C@L{\noexpand\Cal}\def\B@B{\noexpand\Bbb}\def\fR@K{\noexpand\frak}
\def\S@{\noexpand\S}\def\P@P{\noexpand\"}
\def\xpar{\\}

\def\writetoc#1{\iftocloaded\ifThisToToc\begingroup\def\totoc{}
  \def\Cal{\noexpand\C@L}\def\Bbb{\noexpand\B@B}
  \def\frak{\noexpand\fR@K}\def\goth{\frak}\def\S{\noexpand\S@}
  \def\"{\noexpand\P@P}
  \def\xpar{\par\penalty100000 }\def\idx##1{##1}\def\\{\xpar}
  \edef\next@{\write\toc{\noindent#1\leaderfill\noexpand\folio\par}}%
  \next@\endgroup\global\ThisToTocfalse\fi\fi}
\def\leaderfill{\leaders\hbox to 1em{\hss.\hss}\hfill}

\newif\ifindexloaded \indexloadedfalse
\def\idx#1{\ifindexloaded\begingroup\def\ign{}\def\it{}\def\/{}%
 \def\smc{}\def\bf{}\def\tt{}%
 \def\Cal{\noexpand\C@L}\def\Bbb{\noexpand\B@B}%
 \def\frak{\noexpand\fR@K}\def\goth{\frak}\def\S{\noexpand\S@}%
  \def\"{\noexpand\P@P}%
 {\edef\next@{\write\index{#1, \noexpand\folio}}\next@}%
 \endgroup\fi{#1}}
\def\ign#1{}

\def\totoc{\global\ThisToToctrue}

\outer\def\head#1\endhead{\par\penaltyandskip@{-200}\aboveheadskip
 {\headfont@\raggedcenter@\interlinepenalty\@M
 \ignorespaces#1\endgraf}\nobreak
 \vskip\belowheadskip
 \headmark{#1}\writetoc{#1}}

\outer\def\chaphead#1\endchaphead{\par\penaltyandskip@{-200}\aboveheadskip
 {\chapheadfonts\raggedcenter@\interlinepenalty\@M
 \ignorespaces#1\endgraf}\nobreak
 \vskip3\belowheadskip
 \headmark{#1}\writetoc{#1}}

\def\folio{{\foliofont@\ifnum\pageno<\z@ \romannumeral-\pageno
 \else\number\pageno \fi}}
\newtoks\leftheadtoks
\newtoks\rightheadtoks

%Uppercase ist abgestellt:
\def\leftheadtext{\nofrills@{\relax}\lht@
  \DNii@##1{\leftheadtoks\expandafter{\lht@{##1}}%
    \mark{\the\leftheadtoks\noexpand\else\the\rightheadtoks}
    \ifsyntax@\setboxz@h{\def\\{\unskip\space\ignorespaces}%
        \headlinefont@##1}\fi}%
  \FN@\next@}
%Uppercase ist abgestellt:
\def\rightheadtext{\nofrills@{\relax}\rht@
  \DNii@##1{\rightheadtoks\expandafter{\rht@{##1}}%
    \mark{\the\leftheadtoks\noexpand\else\the\rightheadtoks}%
    \ifsyntax@\setboxz@h{\def\\{\unskip\space\ignorespaces}%
        \headlinefont@##1}\fi}%
  \FN@\next@}
%\headline={\def\chapter#1{\chapterno@. }%
%  \def\\{\unskip\space\ignorespaces}\headlinefont@
%  \ifodd\pageno \rightheadline \else \leftheadline\fi}
\def\NoRunningHeads{\global\runheads@false\global\let\headmark\eat@}

\newif\iffirstpage@     \firstpage@true
\newif\ifrunheads@      \runheads@true

%Erg"nzungen zu Runningheads und Pagenumbers:
\newdimen\fullhsize \fullhsize=\hsize
\newdimen\fullvsize \fullvsize=\vsize
\def\fullline{\hbox to\fullhsize}

\def\pagenumbers{\gdef\folio{\folio@}}

\let\norunningheads\NoRunningHeads
\def\userunningheads{\global\runheads@true}
%Default: Seitennumerierung unten.
\norunningheads

\headline={\def\chapter#1{\chapterno@. }%
  \def\\{\unskip\space\ignorespaces}\ifrunheads@\headlinefont@
    \ifodd\pageno\rightheadline \else\leftheadline\fi
   \else\hfil\fi\ifNoRunHeadline\global\NoRunHeadlinefalse\fi}
\let\folio@\folio
\def\foliofont@{\foliofont}
\def\foliofont{\eightrm}
\def\headlinefont@{\headlinefont}
\def\headlinefont{\eightpoint\smc}
\def\leftheadline{\rlap{\folio}\hfill
   \ifNoRunHeadline\else\iftrue\topmark\fi\fi \hfill}
\def\rightheadline{\hfill\ifNoRunHeadline
   \else \expandafter\iffalse\botmark\fi\fi
  \hfill \llap{\folio}}
\footline={{\eightpoint\bottremark}%
   \ifrunheads@\else\hfil{\let\foliofont\tenrm\folio}\fi\hfil}
\def\bottremark{}
 
%Definition von \norunninghead:
\newif\ifNoRunHeadline      
\def\norunninghead{\global\NoRunHeadlinetrue}
\norunninghead

\output={\output@}
%\def\output@{\shipout\vbox{%
% \iffirstpage@ \global\firstpage@false
%  \pagebody \logo@ \makefootline%
% \else \ifrunheads@ \makeheadline \pagebody
%       \else \pagebody \makefootline \fi
% \fi}%
% \advancepageno \ifnum\outputpenalty>-\@MM\else\dosupereject\fi}
%
%Modifizierter Output + Index-Output:
\newif\ifoffset\offsetfalse
\output={\output@}
\def\output@{%
 \ifoffset 
  \ifodd\count0\advance\hoffset by0.5truecm
   \else\advance\hoffset by-0.5truecm\fi\fi
 \shipout\vbox{%
  \makeheadline \pagebody \makefootline }%
 \advancepageno \ifnum\outputpenalty>-\@MM\else\dosupereject\fi}

\def\indexoutput#1{%
  \ifoffset 
   \ifodd\count0\advance\hoffset by0.5truecm
    \else\advance\hoffset by-0.5truecm\fi\fi
  \shipout\vbox{\makeheadline
  \vbox to\fullvsize{\boxmaxdepth\maxdepth%
  \ifvoid\topins\else\unvbox\topins\fi% 
  #1 %
  \ifvoid\footins\else % footnote info is present
    \vskip\skip\footins
    \footnoterule
    \unvbox\footins\fi
  \ifr@ggedbottom \kern-\dimen@ \vfil \fi}%
  \baselineskip2pc
  \makefootline}%
 \global\advance\pageno\@ne
 \ifnum\outputpenalty>-\@MM\else\dosupereject\fi}
 
 \newbox\partialpage \newdimen\halfsize \halfsize=0.5\fullhsize
 \advance\halfsize by-0.5em

 \def\begindoublecolumns{\output={\indexoutput{\unvbox255}}%
   \begingroup \def\line{\fullline}
   \output={\global\setbox\partialpage=\vbox{\unvbox255\bigskip}}\eject
   \output={\doublecolumnout}\hsize=\halfsize \vsize=2\fullvsize}
 \def\enddoublecolumns{\output={\balancecolumns}\eject
  \endgroup \pagegoal=\fullvsize%
  \output={\output@}}
\def\doublecolumnout{\splittopskip=\topskip \splitmaxdepth=\maxdepth
  \dimen@=\fullvsize \advance\dimen@ by-\ht\partialpage
  \setbox0=\vsplit255 to \dimen@ \setbox2=\vsplit255 to \dimen@
  \indexoutput{\pagesofar} \unvbox255 \penalty\outputpenalty}
\def\pagesofar{\unvbox\partialpage
  \wd0=\hsize \wd2=\hsize \hbox to\fullhsize{\box0\hfil\box2}}
\def\balancecolumns{\setbox0=\vbox{\unvbox255} \dimen@=\ht0
  \advance\dimen@ by\topskip \advance\dimen@ by-\baselineskip
  \divide\dimen@ by2 \splittopskip=\topskip
  {\vbadness=10000 \loop \global\setbox3=\copy0
    \global\setbox1=\vsplit3 to\dimen@
    \ifdim\ht3>\dimen@ \global\advance\dimen@ by1pt \repeat}
  \setbox0=\vbox to\dimen@{\unvbox1} \setbox2=\vbox to\dimen@{\unvbox3}
  \pagesofar}

\tenpoint
\catcode`\@=\active

\def\smallheadings{\let\chapheadfonts\tenpoint\let\headfonts\tenpoint}

\tenpoint
\catcode`\@=\active

\def\LL{\leavevmode\setbox0=\hbox{L}\hbox to\wd0{\hss\char'40L}}

            %used for crossreferencing, Tex should ignore.
             %used for refencing (section-numbers)
          %used for new-section numbers

\def\today{\ifcase\month\or
 January\or February\or March\or April\or May\or June\or
 July\or August\or September\or October\or November\or December\fi
 \space\number\day, \number\year}
 %zum Nummerieren

\def\({\left(}
\def\){\right)}
\def\[{\left[}
\def\]{\right]}

\def\3{\ss}
\catcode`\@=11
\def\dddot#1{\vbox{\ialign{##\crcr
      .\hskip-.5pt.\hskip-.5pt.\crcr\noalign{\kern1.5\p@\nointerlineskip}
      $\hfil\displaystyle{#1}\hfil$\crcr}}}

\newif\iftab@\tab@false
\newif\ifvtab@\vtab@false
\def\tab{\bgroup\tab@true\vtab@false\vst@bfalse\Strich@false%
   \def\\{\global\hline@@false%
     \ifhline@\global\hline@false\global\hline@@true\fi\cr}
   \edef\l@{\the\leftskip}\ialign\bgroup\hskip\l@##\hfil&&##\hfil\cr}
\def\endtab{\cr\egroup\egroup}
\def\vtab{\vtop\bgroup\vst@bfalse\vtab@true\tab@true\Strich@false%
   \bgroup\def\\{\cr}\ialign\bgroup&##\hfil\cr}
\def\endvtab{\cr\egroup\egroup\egroup}
\def\stab{\D@cke0.5pt\null 
 \bgroup\tab@true\vtab@false\vst@bfalse\Strich@true\Let@@\vspace@
 \normalbaselines\offinterlineskip
  \openup\spreadmlines@
 \edef\l@{\the\leftskip}\ialign
 \bgroup\hskip\l@##\hfil&&##\hfil\crcr}
\def\endstab{\crcr\egroup
 \egroup}
\newif\ifvst@b\vst@bfalse
\def\vstab{\D@cke0.5pt\null
 \vtop\bgroup\tab@true\vtab@false\vst@btrue\Strich@true\bgroup\Let@@\vspace@
 \normalbaselines\offinterlineskip
  \openup\spreadmlines@\bgroup}
\def\endvstab{\crcr\egroup\egroup
 \egroup\tab@false\Strich@false}

\newdimen\htstrut@
\htstrut@8.5\p@
\newdimen\htStrut@
\htStrut@12\p@
\newdimen\dpstrut@
\dpstrut@3.5\p@
\newdimen\dpStrut@
\dpStrut@3.5\p@
\def\openup{\afterassignment\@penup\dimen@=}
\def\@penup{\advance\lineskip\dimen@
  \advance\baselineskip\dimen@
  \advance\lineskiplimit\dimen@
  \divide\dimen@ by2
  \advance\htstrut@\dimen@
  \advance\htStrut@\dimen@
  \advance\dpstrut@\dimen@
  \advance\dpStrut@\dimen@}
\def\Let@@{\relax\iffalse{\fi%
    \def\\{\global\hline@@false%
     \ifhline@\global\hline@false\global\hline@@true\fi\cr}%
    \iffalse}\fi}
\def\matrix{\null\,\vcenter\bgroup
 \tab@false\vtab@false\vst@bfalse\Strich@false\Let@@\vspace@
 \normalbaselines\openup\spreadmlines@\ialign
 \bgroup\hfil$\m@th##$\hfil&&\quad\hfil$\m@th##$\hfil\crcr
 \Mathstrut@\crcr\noalign{\kern-\baselineskip}}
\def\endmatrix{\crcr\Mathstrut@\crcr\noalign{\kern-\baselineskip}\egroup
 \egroup\,}
\def\smatrix{\D@cke0.5pt\null\,
 \vcenter\bgroup\tab@false\vtab@false\vst@bfalse\Strich@true\Let@@\vspace@
 \normalbaselines\offinterlineskip
  \openup\spreadmlines@\ialign
 \bgroup\hfil$\m@th##$\hfil&&\quad\hfil$\m@th##$\hfil\crcr}
\def\endsmatrix{\crcr\egroup
 \egroup\,\Strich@false}
\newdimen\D@cke
\def\Dicke#1{\global\D@cke#1}
\newtoks\tabs@\tabs@{&}
\newif\ifStrich@\Strich@false
\newif\iff@rst

\def\Stricherr@{\iftab@\ifvtab@\errmessage{\noexpand\s not allowed
     here. Use \noexpand\vstab!}%
  \else\errmessage{\noexpand\s not allowed here. Use \noexpand\stab!}%
  \fi\else\errmessage{\noexpand\s not allowed
     here. Use \noexpand\smatrix!}\fi}
\def\format{\ifvst@b\else\crcr\fi\egroup\iffalse{\fi\ifnum`}=0 \fi\format@}
\def\format@#1\\{\def\preamble@{#1}%
 \def\Str@chfehlt##1{\ifx##1\s\Stricherr@\fi\ifx##1\\\let\Next\relax%
   \else\let\Next\Str@chfehlt\fi\Next}%
 \def\c{\hfil\noexpand\ifhline@@\hbox{\vrule height\htStrut@%
   depth\dpstrut@ width\z@}\noexpand\fi%
   \ifStrich@\hbox{\vrule height\htstrut@ depth\dpstrut@ width\z@}%
   \fi\iftab@\else$\m@th\fi\the\hashtoks@\iftab@\else$\fi\hfil}%
 \def\r{\hfil\noexpand\ifhline@@\hbox{\vrule height\htStrut@%
   depth\dpstrut@ width\z@}\noexpand\fi%
   \ifStrich@\hbox{\vrule height\htstrut@ depth\dpstrut@ width\z@}%
   \fi\iftab@\else$\m@th\fi\the\hashtoks@\iftab@\else$\fi}%
 \def\l{\noexpand\ifhline@@\hbox{\vrule height\htStrut@%
   depth\dpstrut@ width\z@}\noexpand\fi%
   \ifStrich@\hbox{\vrule height\htstrut@ depth\dpstrut@ width\z@}%
   \fi\iftab@\else$\m@th\fi\the\hashtoks@\iftab@\else$\fi\hfil}%
 \def\s{\ifStrich@\ \the\tabs@\vrule width\D@cke\the\hashtoks@%
          \fi\the\tabs@\ }%
 \def\sa{\ifStrich@\vrule width\D@cke\the\hashtoks@%
            \the\tabs@\ %
            \fi}%
 \def\se{\ifStrich@\ \the\tabs@\vrule width\D@cke\the\hashtoks@\fi}%
 \def\cd{\hfil\noexpand\ifhline@@\hbox{\vrule height\htStrut@%
   depth\dpstrut@ width\z@}\noexpand\fi%
   \ifStrich@\hbox{\vrule height\htstrut@ depth\dpstrut@ width\z@}%
   \fi$\dsize\m@th\the\hashtoks@$\hfil}%
 \def\rd{\hfil\noexpand\ifhline@@\hbox{\vrule height\htStrut@%
   depth\dpstrut@ width\z@}\noexpand\fi%
   \ifStrich@\hbox{\vrule height\htstrut@ depth\dpstrut@ width\z@}%
   \fi$\dsize\m@th\the\hashtoks@$}%
 \def\ld{\noexpand\ifhline@@\hbox{\vrule height\htStrut@%
   depth\dpstrut@ width\z@}\noexpand\fi%
   \ifStrich@\hbox{\vrule height\htstrut@ depth\dpstrut@ width\z@}%
   \fi$\dsize\m@th\the\hashtoks@$\hfil}%
 \ifStrich@\else\Str@chfehlt#1\\\fi%
 \setbox\z@\hbox{\xdef\Preamble@{\preamble@}}\ifnum`{=0 \fi\iffalse}\fi
 \ialign\bgroup\span\Preamble@\crcr}
\newif\ifhline@\hline@false
\newif\ifhline@@\hline@@false
\def\hlinefor#1{\multispan@{\strip@#1 }\leaders\hrule height\D@cke\hfill%
    \global\hline@true\ignorespaces}
\def\Item "#1"{\par\noindent\hangindent2\parindent%
  \hangafter1\setbox0\hbox{\rm#1\enspace}\ifdim\wd0>2\parindent%
  \box0\else\hbox to 2\parindent{\rm#1\hfil}\fi\ignorespaces}
\def\ITEM #1"#2"{\par\noindent\hangafter1\hangindent#1%
  \setbox0\hbox{\rm#2\enspace}\ifdim\wd0>#1%
  \box0\else\hbox to 0pt{\rm#2\hss}\hskip#1\fi\ignorespaces}
\def\item"#1"{\par\noindent\hang%
  \setbox0=\hbox{\rm#1\enspace}\ifdim\wd0>\the\parindent%
  \box0\else\hbox to \parindent{\rm#1\hfil}\enspace\fi\ignorespaces}
\let\plainitem@\item
\catcode`\@=13
\hsize13cm
\vsize19cm
\newdimen\fullhsize
\newdimen\fullvsize
\newdimen\halfsize
\fullhsize13cm
\fullvsize19cm
\halfsize=0.5\fullhsize
\advance\halfsize by-0.5em

\magnification1200

\TagsOnRight

\def\AnARAA{1}
\def\CaDLAA{2}
\def\FoWaAA{3}
\def\GaRaAA{4}
\def\KratHYP{5}
\def\SlatAC{6}
\def\WhWaAA{7}

%equation numbers
\def\AA{1}
\def\AB{2}
\def\AC{3}
\def\AD{4}
\def\AE{5}

%theorem numbers

%figure numbers

\def\po#1#2{(#1)_#2}

\topmatter 
\title Asymptotic analysis of a Selberg-type integral via hypergeometrics
\endtitle 
\author C.~Krattenthaler$^{\dagger}$
\endauthor 
\affil 
Fakult\"at f\"ur Mathematik, Universit\"at Wien,\\
Nordbergstra{\ss}e~15, A-1090 Vienna, Austria.\\
WWW: \tt http://www.mat.univie.ac.at/\~{}kratt
\endaffil
\address Fakult\"at f\"ur Mathematik, Universit\"at Wien,
Nordbergstrasze~15, A-1090 Vienna, Austria.
WWW: \tt http://www.mat.univie.ac.at/\~{}kratt
\endaddress
\thanks $^\dagger$Research partially supported 
by the Austrian Science Foundation FWF, grants Z130-N13 and S9607-N13,
the latter in the framework of the National Research Network
``Analytic Combinatorics and Probabilistic Number Theory"%
\endthanks
\subjclass Primary 33C20;
 Secondary 33A15  33C52 60B20 82B05
\endsubjclass
\keywords Selberg-type integral, hypergeometric series, contiguous relation,
balanced hypergeometric series, transformation formula
\endkeywords
\abstract 
We show how to determine the asymptotics of a certain Selberg-type 
integral by means of tools available in the theory of
(generalised) hypergeometric series. This provides an alternative
derivation of a result of Carr\'e, Deneufch\^atel, Luque and Vivo
[{\tt ar$\chi$iv:1003.5996}]. 
\endabstract
\endtopmatter
\document

In \cite{\CaDLAA}, Carr\'e, Deneufch\^atel, Luque and Vivo consider
the Selberg-type integral  
$$
S_k(a,b)\frac {1} {N!}\int _{[0,1]^N} ^{}x_1^k
\Bigg(\prod _{1\le i<j\le N} ^{}(x_i-x_j)^2\Bigg)
\Bigg(\prod _{i=1} ^{N}x_i^{a-1}(1-x_i)^{b-1}\,dx_i\Bigg),
$$
and they
determine its asymptotic behaviour when $N,a,b$ all tend to infinity
so that $a\sim a_1N$ and $b\sim b_1N$, where $a_1$ and $b_1$ are given
non-negative real numbers. The reader is referred to the introduction
of \cite{\CaDLAA} for information on motivation from random matrix
theory connected to random scattering theory to investigate this question.

It should be noted that $S_0(a,b)$ {\it is} a Selberg integral, which can be
evaluated in a product/quotient of gamma functions
(cf.\ \cite{\FoWaAA}). This being the case, the asymptotics of
$S_0(a,b)$ is easily determined by means of known asymptotic formulae
for the Barnes $G$-function (see \cite{\WhWaAA}). 
Thus, it suffices to consider the quotient 
$$J_k=\frac {S_k(a,b)} {S_0(a,b)}$$
(this quotient is denoted by $I_k/N$ in \cite{\CaDLAA}) and determine
its asymptotic behaviour. By (now) classical identities in the theory
of symmetric functions, it is shown in \cite{\CaDLAA, Cor.~II.3} that
$$
J_k=\frac {1} {N\cdot k!}\sum_{i=0} ^{k-1}(-1)^i \binom {k - 1} i
\frac { (N - i)_{ k}\, (a + N - i - 1)_{ k}}
{   (a + b + 2 N - i - 2)_{ k}},
\tag\AA
$$
where $(\alpha)_m$ is the usual Pochhammer symbol defined by
$(\alpha)_m=\alpha(\alpha+1) \cdots(\alpha+m-1)$ for $m\ge1$, and
$(\alpha)_0=1$.
In view of this rather simple expression, which, after substitution of
$a$ by $a_1N$ and of $b$ by $b_1N$, is a rational function in $N$, the
problem of determining the asymptotic behaviour of the resulting
expression looks innocent. 
However, inspection of the expression (\AA) shows that
this is not so innocent at all: the summands are all of the order of magnitude 
$N^k$ for large $N$, while (see Theorem and Corollary below), in fact, the
expression in (\AA) has a {\it finite} limit as $N\to\infty$ in the
described situation. Hence, what happens here is a large amount of 
cancellation of leading asymptotic terms. This is what one has to get under
control, and this constitutes the challenge in determining this limit.

To solve the problem,
Carr\'e, Deneufch\^atel, Luque and Vivo develop a difference
operator calculus in \cite{\CaDLAA}, which is of interest 
in its own right.

The purpose of the present note is to demonstrate that the
tools for solving this problem
are readily available in the literature on
({\it generalised}) {\it hypergeometric series},
and it is clear that they may serve as well to solve problems of
similar type. 

Namely, using the standard hypergeometric notation
$$
{}_p F_q\!\left[\matrix a_1,\dots,a_p\\
b_1,\dots,b_q\endmatrix; z\right]=\sum _{m=0} ^{\infty}\frac
{\po{a_1}{m}\cdots\po{a_p}{m}} {m!\,\po{b_1}{m}\cdots\po{b_q}{m}}
z^m \,
$$
where $\po{a_1}{m}$, etc.\ are again Pochhammer symbols,
Equation~(\AA) becomes
$$
J_k=
\frac {(N+1)_{ k-1}\,
  (a + N-1)_{ k}} {k!\,(2N+ a + b -2)_{ k}}
{}_4 F_3\!\left[\matrix 1 - N, 1 - k, 2 - a - N, 3 - a - b - k - 2N\\
   2 - a - k - N, 1 - k - N, 3 - a - b - 2N\endmatrix; 1\right].
\tag\AB$$
The point of this conversion is that one may now look in the
literature of hypergeometric series, which provides a large arsenal of
utilities, such as transformation formulae for hypergeometric series
etc., and see whether one (or a combination of more) of these formulae
enable us to convert (\AB) into a different form from which the
asymptotics of $J_k$ can be directly read off. As I show below, this
is indeed the case. To the informed reader it will come as no surprise that the
computations below have been carried out using the {\sl Mathematica}
package HYP \cite{\KratHYP}. 

The crucial observation to make is that the $_4F_3$-series in (\AB) 
is not just an ``arbitrary"
hypergeometric series, but that it is a {\it $(-1)$-balanced\/} series
(cf.\ \cite{\AnARAA, p.~140}), that is, that the sum of the lower
parameters of the hypergeometric series minus the sum of
the upper parameters is $-1$. Unfortunately, there are no transformation
formulae available for $(-1)$-balanced series, but there are plenty of
them available for $1$-balanced series (or, as is common to say,
for {\it
balanced series}), such as, for example
(see \cite{\GaRaAA, Appendix (III.16), $q\uparrow1$}; this is T4302 in HYP),
$$\multline
{} _{4} F _{3} \!\left [ \matrix { A, B, C, -n}\\ { E, F, 1 + A + B + C - E -
   F - n}\endmatrix ; {\displaystyle 1}\right ] \\=
{\frac {({ \textstyle A)_n\,(  E + F-A - B)_n\, (  E + F-A - C}) _{n}}
     {({ \textstyle E)_n\, ( F)_n\, ( E + F-A - B - C}) _{n}}}\\
\times
  {} _{4} F _{3} \!\left [ \matrix { -n, E-A,  F-A,  E + F-A - B - C}\\
    { E + F-A - B,  E + F-A - C, 1 - A - n}\endmatrix ; {\displaystyle
    1}\right ]   ,
\endmultline\tag\AC$$
where $n$ is a non-negative integer.
That is, we are just by ``$2$ off" in applying a transformation formula.

In order to relate series with (small) integral differences between
the parameters, the hypergeometric literature offers {\it contiguous
relations}. One of the simplest such relations is (this is C14 in HYP),
$$
  {} _{4} F _{3} \!\left [ \matrix A, B, C, D\\E, F, G
\endmatrix ; {\displaystyle
    z}\right ]  
= z\frac {BCD} {EFG}
  {} _{4} F _{3} \!\left [ \matrix A, B+1, C+1, D+1\\E+1, F+1, G+1
\endmatrix ; {\displaystyle z}\right ] + 
  {} _{4} F _{3} \!\left [ \matrix  A-1, B, C, D\\E, F, G\endmatrix ; {\displaystyle
    z}\right ]   .
$$
If we iterate this contiguous relation, then we arrive at
$$\multline 
  {} _{4} F _{3} \!\left [ \matrix A,B,C,D\\E,F,G
\endmatrix ; {\displaystyle z}\right ] =
        z^r \frac {(B)_r\, (C)_r\, (D)_r} {(E)_r\, (F)_r\, (G)_r} 
  {} _{4} F _{3} \!\left [ \matrix 
A,B+r,C+r,D+r\\E+r,F+r,G+r\endmatrix ; {\displaystyle z}\right ]\\ +
\sum_{s=0} ^{r-1} z^s \frac {(B)_s\, (C)_s\, (D)_s} 
{(E)_s\, (F)_s\, (G)_s} 
   {} _{4} F _{3} \!\left [ \matrix 
A-1,B+s,C+s,D+s\\E+s,F+s,G+s
\endmatrix ; {\displaystyle z}\right ].
\endmultline
\tag\AD$$
We now apply (\AD) with $r=k-1$ to the $_4F_3$-series in (\AB). 
The reader should observe that thereby all the $_4F_3$-series in the
sum over $s$ will be $0$-balanced, while the single $_4F_3$-series
outside of the sum will reduce to $1$ because one of the upper
parameters is $1-k+r=0$. To be precise, 
after writing the result in compact form,
application of (\AD) to (\AB) yields
$$\align 
J_k&=
 \sum _ {s=0} ^{k-1}
\frac {(N+1)_{k-1}\, (a + N-1)_k\, 
      (1 - k)_s\, (2 - a - N)_s\, (3 - a - b - k - 2 N)_s}
   {k!\,( a
   + b + 2 N-2)_k\,(2 - a - k - N)_s\, (1 - k - N)_s\, (3 - a - b - 2 N)_s}\\
&\kern40pt
\times
  {} _{4} F _{3} \!\left [ \matrix 
-N, 1 - k + s, 2 - a - N + s, 3 - a - b - k - 2 N + s\\
       2 - a - k - N + s, 1 - k - N + s, 3 - a - b - 2 N + s
\endmatrix ; {\displaystyle 1}\right ]
\endalign$$
Continuing in the same spirit, if we apply the contiguous relation
(\AD) a second time, this time with $r=k-s-1$, then this will produce
a double sum in which all appearing $_4F_3$-series are balanced. To be
precise, we obtain
$$\multline 
J_k\\=
\sum _{s=0} ^{k-1} \sum _{t=0} ^{k-s-1}
\frac { (N+1)_{k-1}\, (a + N-1)_k\,
(1 - k)_{s+t}\,   (2 - a - N)_{s+t}\, (3 - a - b - k - 2 N)_{s+t}}
 {k!\,(a + b + 2 N-2)_k\,
  (2 - a - k - N)_{s+t}\,(1 - k - N)_{s+t}\,(3 - a - b - 2 N)_{s+t}  
    }\\
\times
  {} _{4} F _{3} \!\left [ \matrix 
3 - a - b - k - 2 N + s + t, -1 - N, 2 - a - N + s + t, 
      1 - k + s + t\\3 - a - b - 2 N + s + t, 1 - k - N + s + t, 
      2 - a - k - N + s + t\endmatrix ; {\displaystyle 1}\right ].
\endmultline$$
Now the transformation formula (\AC) can be applied, and, after some
simplification, the resulting expression collapses to
$$\align
J_k&=
\sum _{s=0} ^{k-1} \sum _{t=0} ^{k-s-1}
\frac { (a-1)_{ k - s - t-1}\, (1 - a - N)_{ s + t+1}\, 
  (k - s - t)_{ s + t}\, (s + t+2)_{ k - s - t-1}} 
{k!\,(2 - a - b - 2 N)_{ k}}\\
&\kern3cm
\cdot
  {} _{4} F _{3} \!\left [ \matrix 
1 - k + s + t, k, a + b + N-2, a + N\\s + t+2, a-1, 
    a + b + 2 N-1\endmatrix ; {\displaystyle 1}\right ]\\
&=\sum _{m=0} ^{k-1} 
\frac { (a-1)_{ k - m-1}\, (1 - a - N)_{ m+1}\, 
  (k - m)_{ m}\, (m+1)_{ k - m}} 
{k!\,(2 - a - b - 2 N)_{ k}}\\
&\kern3cm
\cdot
  {} _{4} F _{3} \!\left [ \matrix 
- k + m+1, k, a + b + N-2, a + N\\m+2, a-1, 
    a + b + 2 N-1\endmatrix ; {\displaystyle 1}\right ].
\endalign$$
Now the limit $N,a,b\to\infty$, in which $a\sim a_1N$ and $b\sim
b_1N$, can be safely performed. 
Hence, we have proven the following result.

\proclaim{Theorem}
The limit of the quantity $J_k$ as $N,a,b\to\infty$ such that $a\sim
a_1N$ and $b\sim b_1N$ is equal to
$$\multline 
\lim_{N\to\infty}J_k=
\sum _{m=0} ^{k-1} 
(-1)^{k-m-1}\binom {k-1}m
\left(\frac {a_1} {a_1+b_1+2}\right)^k
\left(\frac {a_1+1} {a_1}\right)^{m+1}\\
\cdot
\sum _{\ell=0} ^{k-m-1}
(-1)^\ell\binom {k-m-1}\ell
\frac {(k+\ell-1)!(m+1)!} {(k-1)!(m+\ell+1)!}
\left(\frac {(a_1+1)(a_1+b_1+1)} {a_1(a_1+b_1+2)}\right)^\ell.
\endmultline
\tag\AE$$
\endproclaim

This theorem does not give the same expression as in \cite{\CaDLAA,
Theorem~IV.4}. However, once one has found {\it some} expression, it
is easy to pass to an equivalent one by using hypergeometric
transformation and summation formulae. (Since, in the above derivation,
there is plenty of room for variation --- for example, the transformation
formula (\AC) could be applied in different ways, or we could have applied
a different transformation formula for balanced $_4F_3$-series, see
\cite{\SlatAC, Sec.~4.3.5} ---, it is conceivable that there is a direct 
route leading to the corollary below, using the above ideas. We did not
pursue this, however; in particular, we did not do a systematic search
of such an alternative route.)

\proclaim{Corollary \smc(\cite{\CaDLAA, \smc Theorem~IV.4})}
The limit of the quantity $J_k$ as $N,a,b\to\infty$ such that $a\sim
a_1N$ and $b\sim b_1N$ is equal to
$$
\lim_{N\to\infty}J_k=
\frac {1} {k}\sum _{j=0} ^{k-1}(-1)^j 
\binom {k+j-1}j
\frac {(a_1+1)^{j+1}} {(a_1+b_1+2)^{k+j}}
\sum_{i=0} ^{k-j-1}\binom ki\binom k{i+j+1}
(a_1+1)^i.
$$

\endproclaim

\demo{Proof}
In the expression on the right-hand side of (\AE) we interchange the
summations over $m$ and $\ell$, and we write the (now) inner sum over $m$ in
hypergeometric notation. This leads to
$$\multline 
\lim_{N\to\infty}J_k=
\sum _{\ell=0} ^{k-1}(-1)^{k-\ell-1}
\frac {a_1^{k-\ell-1}\,(a_1+1)^{\ell+1}\,(a_1+b_1+1)^{\ell}}
{(a_1+b_1+2)^{k+\ell}}
    \frac {(k+\ell-1)!} {
    \ell !\,(\ell +1)!\,(k-\ell-1)!}\\
\cdot
  {} _{2} F _{1} \!\left [ \matrix 
2, - k + \ell +1\\\ell +2\endmatrix ; {\displaystyle  \frac
{a_1+1} {a_1} }\right ].
\endmultline$$
To the $_2F_1$-series we apply the transformation formula
(see \cite{\SlatAC, (1.7.1.3), sum
reversed at the right-hand side}; this is T2106 in HYP)
$$
{} _{2} F _{1} \!\left [ \matrix { A, -n}\\ { C}\endmatrix ; {\displaystyle
   z}\right ] = {z^n}
  {{  ({ \textstyle C-A}) _{n} }\over {({ \textstyle C}) _{n} }}
  {} _{2} F _{1} \!\left [ \matrix { -n, 1 - C - n}\\ {
       1 + A - C - n}\endmatrix ; {\displaystyle -{{1 - z}\over z}}\right ]  ,
$$
where $n$ is a non-negative integer.
Writing the resulting $_2F_1$-series as a sum over $i$, this gives
$$\multline 
\lim_{N\to\infty}J_k=
\sum _{\ell=0} ^{k-1}
\sum _{i=0} ^{k-\ell-1}
(-1)^{k-\ell-i-1}
\binom {k-\ell-1}i
\frac {(k+\ell-1)!  }
{(\ell-1)!\,\ell!\,(k-\ell-1)!\,(k-i-1)(k-i)
   }\\
\cdot
\frac {(a_1+1)^{k-i}\,(a_1+b_1+1)^{\ell}}
{(a_1+b_1+2)^{k+\ell}}.
\endmultline$$
Now we replace the term $(a_1+b_1+1)^{\ell}$ by
$\sum _{j=0} ^{\ell}\binom \ell {j}(-1)^j(a_1+b_1+2)^{\ell-j}$, we
make the sum over $\ell$ the inner-most sum, and finally write the sum
over $\ell$ in hypergeometric notation. In this way, the last
equation becomes
$$\multline 
\lim_{N\to\infty}J_k=
\sum _{j=0} ^{k-1}
\sum _{i=0} ^{k-j-1}
(-1)^{k-i-1}
\binom {k-j-1}i
\frac {(k+j-1)!  }
{(j-1)!\,j!\,(k-j-1)!\,(k-i-1)(k-i)
   }\\
\cdot
\frac {(a_1+1)^{k-i}}
{(a_1+b_1+2)^{k+j}}
  {} _{2} F _{1} \!\left [ \matrix {-k+i+j+1,k+j }\\ {
j       }\endmatrix ; {\displaystyle 1}\right ]  .
\endmultline$$
The $_2F_1$-series can be evaluated by the Chu-Vandermonde summation
(see \cite{\SlatAC, (1.7.7); Appendix (III.4)}; this is S2101 in HYP)
$$
{} _{2} F _{1} \!\left [ \matrix { A, -n}\\ { C}\endmatrix ; {\displaystyle
   1}\right ] = {{({ \textstyle C-A}) _{n} }\over 
    {({ \textstyle C}) _{n} }},
$$
where $n$ is a non-negative integer. After finally replacing $i$ by
$k-i-j-1$, the claimed result follows immediately upon minor
rearrangement of terms.\quad \quad \qed
\enddemo

\enddocument